\documentclass[11pt, leqno]{article}

\usepackage{amsfonts}

\usepackage{amsmath}

\usepackage{amssymb}

\usepackage{theorem}

{\theorembodyfont{\rmfamily} \newtheorem{Def}{Definition}[section]

\newtheorem{Theorem}{Theorem}[section]

\newtheorem{Lemma}{Lemma}[section]

\newtheorem{Remark}{Remark}[section]

}

\begin{document}

\title{\bf\large Metric Trees, Hyperconvex Hulls and Extensions} \author{\bf A. G. Aksoy \\
\small Department of Mathematics \\ \small Claremont McKenna College \\ \small Claremont, CA 91711 \\
\small asuman.aksoy@claremontmckenna.edu \and \bf B.  Maurizi \\
\small Department of Mathematics \\ \small Washington University in
St. Louis
\\ \small St. Louis, MO 63130. \\ \small bmaurizi@math.wustl.edu }

\date{\small October 22, 2005}

\maketitle



\begin{abstract}

In this paper we examine the relationship between hyperconvex hulls
and metric trees. After providing a linking construction for
hyperconvex spaces, we show that the four-point property is
inherited by the hyperconvex hull, which leads to the theorem that
every complete metric tree is hyperconvex. We also consider some
extension theorems for these spaces.

\end{abstract}



{\bf Keywords:}

Hyperconvex spaces, metric trees, extensions



{\bf AMS subject classification:}

05C12, 54H12, 46M10



\section{Introduction}
    The purpose of this paper is to clarify the relationship between
 metric trees and hyperconvex metric spaces. We provide  a new so-called
 \emph{linking construction} of hyperconvex spaces and show that the four-point
 property of a metric space is inherited by the hyperconvex hull of that
 space. We prove that all complete metric trees are hyperconvex. This in
 turn suggests a new approach to the study of extensions of
 operators. For a metric space $(X,d)$ we use $B(x;r)$ to denote the
 \emph{closed} ball centered at $x$ with radius $r\geq0.$

 \begin{Def}\label{definition}
 A metric space $(X,d)$ is said to be \emph{hyperconvex} if
 $\ \bigcap_{i\epsilon I } B(x_{i};r_{i})\neq\phi $  for every collection $B(x_{i};r_{i})$ of closed
 balls in $X$ for which $d(x_{i},x_{j})\leq r_{i}+r_{j} $.
\end{Def}
 This notion was first introduced by Aronszajn and Panitchpakdi
 in \cite{ref1}, where it is shown that a metric space is hyperconvex
 if and only if it is injective with respect to nonexpansive
 mappings. Later Isbell \cite{ref7} showed that every metric space has
 an injective hull, therefore every metric space is isometric to a
 subspace of a minimal hyperconvex space. Hyperconvex metric spaces are complete and connected
 \cite{ref9}.
  The simplest examples of
 hyperconvex spaces are the set of real numbers $\mathbb R$, or a
 finite-dimensional real Banach space endowed with the maximum norm.
 While the Hilbert space $l^2$ fails to be hyperconvex, the
 spaces $L^\infty$ and $l^\infty$ are hyperconvex. In \cite{ref2} it
 is shown that $\mathbb{R}^2$ with the ``river" or ``radial" metric is
 hyperconvex. We will show that there is a general ``linking
 construction" yielding hyperconvex spaces. Constructions of the river and radial
 metrics are obtained as  special cases. Moreover, in these spaces
 paths between points are restricted; they must pass through
 certain ``common" points. On the other hand, the concept of a metric
 tree in graph theory also has a built-in restriction. A complete
 metric space $X$ is a \emph{metric tree} provided that for any two points $x$
 and $y$ in $X$ there is a unique arc joining $x$ and $y$, and this
 arc is a geodesic arc. For more on metric trees we refer the reader to \cite{ref3}, \cite{ref5},\cite{ref6} and
 \cite{ref13}.
 One particularly useful characterization of
 metric trees is given by the ``four-point condition".
\begin{Def}\label{definition}
 A metric space $(X,d)$ is said to satisfy the  \emph{four-point} property provided that for each set of four points
 $x,y,u,v$ in $X$ the following holds:

$$
 d(x,y)+d(u,v)\leq \max (d(x,u)+d(y,v),d(x,v)+d(y,u)) .
$$

\end{Def}

The four-point condition is stronger than the triangle inequality
(take $u=v$), but it should not be confused with the ultrametric
definition. An ultrametric satisfies the condition $d(x,y)\leq
\max(d(x,z),d(y,z))$, and this is stronger than the four-point
condition. The four-point condition is equivalent to saying two of
the three numbers
$$d(x,y)+d(u,v),\qquad d(x,u)+d(y,v),\qquad d(x,v)+d(y,u)$$ are the same and the
third one is less than or equal to that number.\, The study of
spaces with the four-point property has a practical motivation (in
numeric taxonomy), but also has interesting theoretical aspects. If
the space $X$ is finite then $X$ can be imagined as subspaces of
usual graph-theoretic trees (with nonnegative weight on edges
determining their length). In \cite{ref5} it is shown that a metric
space is a metric tree if and only if it is complete, connected and
satisfies the four-point property. The first section of this paper
is devoted to hyperconvex spaces and hyperconvex hulls. Next we show
that the four-point property is inherited by the hyperconvex hull.
In the last section, we mention some known extension properties in
the context of $P_1$-spaces, which can be rephrased now for complete
metric trees.

\section{The Linking Construction for Hyperconvex Spaces, and
the Hyperconvex Hull}

    The understanding of hyperconvex spaces rests on how these
spaces can be constructed. There is one obvious way to construct a
hyperconvex space which is analogous to the direct product: take a
collection of hyperconvex spaces and put the supremum metric on the
Cartesian product. This new space will be hyperconvex essentially
because any pairwise overlapping collection must overlap in each
coordinate. In the following we will present two different
constructions, each of which builds a larger space out of smaller
spaces.  We will take several hyperconvex spaces and join each of
them by one point to a central hyperconvex space. This type of
linking creates a restrictive  movement in the sense that in order
to pass between different points in different spaces, one must
travel through the common point, and through the central hyperconvex
space. A similar construction to this is also presented in
\cite{ref2} and \cite{ref8}.

Consider a metric space $(X,d)$ and an arbitrary set $C$ outside the
set $X$. Let $f:X \cup C \rightarrow X\times [0,\infty)$  be such
that $f_1=f_{|X}=(x,0)$ and $f_2=f_{|C}:C\rightarrow X\times
(0,\infty)$. The first coordinate can be thought of as the closest
point in $X$ to the point in the domain, and the second coordinate
can be thought of as the distance to that closest point. Let us
define a metric $\rho$ on $X \cup C $ as follows:
\begin{displaymath}
\rho(p_{1}, p_{2})=\left\{ \begin{array}{ll}
                    0, & \mbox{if $p_{1}=p_{2}$};\\[.04in]
                    d(p_{1},p_{2}), & \textrm{if $p_{1},p_{2}
                    \in X$};\\[.04in]
                    d(f_1(p_1),p_2)+f_2(p_1), & \textrm{if $p_1\in C,\, p_2\in
                    X$};\\[.04in]
                    d(f_1(p_1),f_1(p_2))+f_2(p_1)+f_2(p_2)
                    &\textrm{if $p_1,\, p_2\in C$}.
                    \end{array} \right.
\end{displaymath}
It is straightforward to check that $\rho$ is a metric.\\

In the following theorem we think of $f(\alpha)$ as the point in $X$
at which the entire space $W_{\alpha}$ is linked to $X$, and
$g(\alpha)$ as the point in $W_{\alpha}$ at the other end of that
link.
\begin{Theorem}\label{Theorem}
Suppose $(X,d)$ is a hyperconvex metric space and ${(W_\alpha\,
d_\alpha)}_{\alpha\in I}$ is a collection of hyperconvex spaces.
Given a function $f:I\rightarrow X $ and a function $g:I\rightarrow
W_{\alpha}\backslash \{g(\alpha)\} $ such that $g(\alpha)\in
W_\alpha $, one can construct a metric $\rho(x,y)$ where
\begin{displaymath}
\rho(x,y)=\left\{\begin{array}{ll}
          d(y,z), & \textrm{for $y,z,\in X$}; \\
          d(y,f(\alpha))+d_{\alpha} (g(\alpha),z), & \textrm{for $y\in X, z\in W_{\alpha}\backslash g(\alpha)$};\\
          d(f(\beta),f(\alpha))+d_{\alpha}(g(\alpha),y)+d_{\beta}(g(\beta),z), &
          \textrm{for $y\in W_{\alpha}\backslash g(\alpha), z\in
          W_{\beta}\backslash g(\beta)$}.
          \end{array}\right.
\end{displaymath}
is a metric on the set $Z:=X\cup (W_{\alpha}\backslash
\{g(\alpha)\})$ such that it is hyperconvex.
\end{Theorem}
\textbf{Proof:}  Consider a hyperconvex metric space $(X,d)$, and a
set $C= \cup (W_{\alpha}\backslash g\{(\alpha)\})$. Use the above
construction to define a function $$F:X\cup ({W_{\alpha}\backslash
\{g(\alpha)\}}) \rightarrow X\times [0,\infty)$$
\\by $F(x)=(x,0)$ for $x\in X$ \,and\,
$F(y)=(f(\alpha),d_{\alpha}(y,g(\alpha))$ for $y\in \cup
(W_{\alpha}\backslash \{g(\alpha)\}).$ \\ [.05in] Notice that
$d(y,g(\alpha))>0$ for all $y$. Therefore, the metric on $Z:={X\cup
(W_{\alpha}\setminus \{g(\alpha)\})}$ is exactly as the one stated.
To prove $Z$ is hyperconvex we consider two cases. In the first case
we assume balls ``overflow" into $X$ which is hyperconvex; in the
second case one of the balls does not overflow into $X$ so the total
intersection must be found in $W_{\alpha}\backslash
\{g(\alpha)\}$.\\[.05in]
Case 1: Let $r_{i}\geq \rho(x_i,F_1(x_i))$ for all $i$, and let $
\overline{r_j}= r_j-\rho(x_j,F_1(x_j))$ . Now notice that
$$\rho(x_i,x_k)=\rho(x_j,F_1(x_j))
+\rho(F_1(x_j),F_1(x_k))+\rho(F_1(x_k),x_k)\leq r_j+r_k .$$ This
implies
$\rho(F_1(x_j),F_1(x_k))\leq\overline{r_j}+\overline{r_k}$.\, Since
$X$ is hyperconvex and $F_1(x_j),F_1(x_k)\in X$ we have
$\bigcap_{i\in I} B(F_1(x_i),\overline {r_i}) \neq \phi.$ However we
already have $B(F_1(x_i),\overline{r_i})\subset B(x_i,r_i)$.\\
[.07in] Case 2: Suppose we have $x_m\in W_\alpha \setminus
\{g(\alpha)\}$ with $r_m$ such that $r_m < \rho {(x_m,F_1(x_m))} $.
Now observe that for any $x_i\not\in W_{\alpha} \setminus
\{g(\alpha)\}$, we have
$$\rho (x_m,x_i) = \rho (x_m,F_1(x_m)) +\rho(F_1(x_m),x_i) \leq
r_m+r_i,$$ and this together with the condition on $r_m$ implies
that $ \rho (F_1(x_m),x_i) < r_i$. We now set $$\overline{r_I} =
r_i-\rho(F_1(x_m),x_i) ~\mbox{and}~  J:=\{i\in I: x_i\not\in
W_{\alpha} \setminus g(\alpha)\}.$$ Since $r_j > 0$  we have $
\bigcap_{j\in J} B(g(\alpha) , \overline r_j) \neq\phi$ and from
hyperconvexity of $W_\alpha$ we also know  $\bigcap_{i\in I\setminus
J}B(x_i,r_i)\neq\phi$. Note that $g(\alpha)\not\in B(x_m,r_m),$
therefore the intersection point cannot be $g(\alpha)$.  Next we
claim that balls of the form  $B(g(\alpha), \overline {r_j}) $ where
$j\in J$, and $B(x_i,r_i)$ where $i\in I\setminus J$, will intersect
pairwise. For if we consider
$$\rho(x_j,F_1(x_m))+d_\alpha (g(\alpha),x_i)= \rho(x_i,x_j)\leq r_i
+ r_j,$$ subtracting $ \rho(x_j,F_1(x_m) )$ from both sides will
give  $$d_\alpha (g(\alpha),x_i) \leq r_i+ r_j -\rho(x_i,x_j)=
r_i+\overline r_j .$$ Using the hyperconvexity of $W_\alpha$,
$$\left[\cap_{j\in J} B(g(\alpha),\overline r_j)\right] \cap \left[\cap_{i\in
I\setminus J} B(x_i,r_j)\right]\neq\phi .$$ Finally, noting
$B(g(\alpha), \overline r_j) \setminus \{g(\alpha)\} \subset
B(x_j,r_j)$, we have $ \cap_{i\in I} B(x_i,r_i)\neq\phi$. This
concludes the proof.

\hfill $\Box$

Next we show a way to construct a hyperconvex space from a given
normed space by defining a different metric on this space. We take
an appropriate subspace having a hyperconvex metric, and then
decompose the normed space into subspaces linked with all rays
connecting points outside the subspace with their closest point.
First, we need the following lemma which illustrates that if we have
a subspace of a normed space for which the closest point exists and
is unique, then one can partition the remaining points of the space
into equivalence classes, by defining two points to be equivalent if
they lie on the same ray from the subspace.

\begin{Lemma}
Suppose $X$ is a normed space and $Z$ is a subspace such that the
closest point in $Z$ to any $x\in X$ exists and is unique. Suppose
$h(p)$ is the closest point in $Z$ to $p\in X\setminus Z$, and the
ray pointing from $h(p)$ in the direction of $p$ is denoted by
$\lambda_p$ (i.e., $ \lambda_p= \mu (p-h(p))+h(p))$ where $\mu \in
[0,\infty)$). Then, if $p\in \lambda_q$,  $p=t_0(q-h(q))+h(p)$
implies $h(p)=h(q)$.
\end{Lemma}
\textbf{Proof}: Suppose $ t_0 <1 $. Let $ z\in Z $, $ z\neq h(q)$.
Then we have
$$
d(h(q),q)=d(h(q),p)+d(p,q)\leq d(z,p)+d(p,q)
 $$
 where the first equality comes
from the fact that $p$ lies on a line segment between $q$ and
$h(q)$, and the second inequality is a consequence of the fact that
$h(q)$ is minimal and unique. Therefore, we have $d(h(q),q) <
d(z,q).$ For the case $t_0>1$, suppose that for some $z\in Z$,
$d(z,p)<d(h(q),p)$.
\\Let $$\beta=\frac {d(q,h(q))}{d(p,h(q))}.$$  Set \, $ z^\ast
=(1-\beta)h(q)+\beta z$, and compute $$ d(z^\ast,q) =  \parallel
(1-\beta)h(q)+\beta z-[h(q)+\beta(p-h(q))]\parallel =  \parallel
\beta (z-p)\parallel = \beta\, d(z,p).$$ This yields  $ d(
z^\ast,q)<d(h(q),q)$, giving a contradiction. This means that $p\in
\lambda _q $ implies $\lambda _p=\lambda _q .$ \hfill $\Box$

\begin{Theorem} Suppose $X$ is a normed space and $Z$ is a subspace
such that the closest point to any $X$ exists and is unique. Suppose
also that $Z$ has a different metric with which $Z$ is hyperconvex.
One can construct a metric on $X$ so that it is hyperconvex.
\end{Theorem}
\textbf{Proof:}  Consider the equivalence classes of rays
$[{\lambda_\alpha}]_{\alpha \in I}$ described in the above lemma. We
have the functions $f:I\rightarrow Z$  which takes $\alpha \mapsto
h(p)$ for a point $p\in \lambda_\alpha $ and $g: I\rightarrow
\bigcup_{\alpha\in I} \lambda_\alpha $ which takes $\alpha \mapsto
h(p)$. We assumed that $Z$ is hyperconvex under some metric
$\delta$. Each of the $\lambda_\alpha$ is a hyperconvex metric space
under the norm restricted to $\lambda _\alpha$, since  $\lambda
_\alpha$ is isometric to $ [0,\infty)$. By  Theorem 2.1 we have a
hyperconvex space {$ Z \bigcup\, \cup_{\alpha \in I} (\lambda_\alpha
\setminus p_\alpha)$}. However, this the normed space $X$, with the
metric

\begin{displaymath}
d(x,y)=\left\{\begin{array}{ll}
       \parallel x-y\parallel,  & \textrm{if $ x,y
       \in \lambda_
       \alpha \setminus h(x)$};\\[.04in]
       \delta(x,y),

        & \textrm{if $x,y \in Z$}; \\[.04in]
       \delta(x,h(y))+\parallel y-h(y)\parallel, & \textrm {if $ x \in
       Z, y \in \lambda_\alpha \setminus h(y)$}; \\[.04in]
       \delta(h(x),h(y)) + \parallel x-h(x) \parallel+\parallel
       y-h(y) \parallel, & \textrm {if $x\in \lambda_\alpha \setminus
     h(x), y\in \lambda_\beta \setminus h(y)$}.
     \end{array}\right.
\end{displaymath}
\hfill $\Box$

Notice that if $ X= \mathbb {R}^2$ and $ Z$ is the $x$-axis, then
this  metric is the``river metric", and if $ X= \mathbb{R}^2$ and $
Z=(0,0)$
then it is the ``radial metric" described in \cite{ref2}. \\[.07in]
\begin{Def} \label{definition}
 Given a metric space $(X,d)$, the \emph{hyperconvex hull} of $X$ is
another metric space $ (Y,\rho)$ such that $X$ is contained
isometrically in $Y$, where $Y$ is a hyperconvex metric space and Y
is minimal.
\end{Def}

It is not immediately clear that such a metric space exists or is
unique. Given a collection of points $\{ x_i\}_{i\in I} \in X$ and
radii $\{r_i\}_{i\in I} \in \mathbb{R}^+$, we say that this
collection is \emph{pairwise overlapping} if $$d(x_i,x_j)\leq
r_i+r_j$$ for all $i,j \in I$. In a given metric space $(X,d)$, if
we have an overlapping collection $\{x_\alpha \}_{\alpha \in I} \in
X $, $\{r_\alpha\}_{\alpha \in I} \in \mathbb {R}^+$ we can shrink
any overlapping collection until it is minimal. We say it is
\emph{minimally overlapping} if for all $\epsilon
>0$ and for all $\beta\in I$, the collection of points
$$\{x_\alpha\}_{\alpha \in I}    ,\qquad \{r_\alpha\} _{\alpha \in I,
\alpha\neq \beta} \cup [r_\beta -\epsilon ]$$ is not pairwise
overlapping. In other words, minimally overlapping means we can not
shrink any of the radii. \,Now using a Zorn's lemma argument, for
any pairwise overlapping collection $\{x_\alpha\}_{\alpha\in I},
\{r_\alpha\}_{\alpha\in I}$ with $x_\alpha\in X$, $r_\alpha \in
\mathbb {R}^+$, we can find a set of radii $\{r^*_\alpha\}_{\alpha
\in I}$  with  $r^*_\alpha\leq r_\alpha$  such that the collection
$\{x_\alpha\}$, $\{r^*_\alpha\}_{\alpha \in I}$ is minimally
overlapping. \,Analogous to the completion of a metric space, to
construct a hyperconvex
 hull one takes  a pairwise overlapping collection with no total
 intersection, and regards it as a single object in the set of all
 such objects. Then, putting a suitable metric on this set results in
 a metric space with the desired property. In the following we will
denote the hyperconvex hull by $h(X)$.
\begin{Def}\label {Definition}
A function $f\in C(X)$ is called a \emph{minimal extremal function }
if $$f(x)+f(y)\geq d(x,y),$$ and is pointwise minimal. That is, if
$g$ is another function with the same property such that $g(x)\leq
f(x)$ for all $x\in X$, then $g=f$.
\end{Def}

The similarity between a minimally overlapping collection and a
minimal extremal function is explained in the following remark.
\begin{Remark} Suppose we have a minimally overlapping
collection $\{ x_\alpha\} \in X$ and $\{r_\alpha\}\in \mathbb {R}^+
$. We can think of this collection as a function
$$\tilde{f}:\{x_\alpha\}_{\alpha \in I} \rightarrow \mathbb {R}^+$$
defined by $ x_\alpha\mapsto r_\alpha $. Because of a pairwise
overlap we have $\tilde{f}(x)+\tilde{f}(y) \geq d(x,y)$. Moreover,
we can extend $\tilde{f}$ to $f$ where $$f:X\rightarrow \mathbb
{R}^+$$ and
$$f(x)+f(y) \geq d(x,y).$$ To do this, we define
$f:X\rightarrow \mathbb{R}^+$ by $$x\mapsto ~\mbox{inf}~ _{x_\alpha}
[d(x,x_\alpha)+r_\alpha].$$ It is easy to show that $f$ is extremal \cite{ref7}.\\

There is an obvious family of minimal extremal functions on $X$,
namely, select $x\in X$ and define a function $h_x$ by:
$$h_x(z)=d(x,z)$$ obviously $h_x(x)=0$. We will call these
\emph{distance cones}. One natural question is  whether or not there
are other minimal extremal functions besides distance cones? The
answer to this question is in the connection between hyperconvexity
and minimal extremal functions. It was shown by Isbell \cite{ref7}
that there are other extremal minimal functions precisely when the
space is \underline{not} hyperconvex. The following theorem (proof
can be found in \cite{ref7}) introduces the basic properties of the
hyperconvex hull.
\end{Remark}

\begin{Theorem}
For a metric space $(X,d)$, consider the set
$$h(X)= \{ f:X\rightarrow \mathbb{R}:f(x)+f(y)\geq d(x,y) \textrm{
and $f$ is minimal}\}$$ and the metric $$\rho(f,g)=\sup_{x\in X}
{d(f(x),g(x))}$$ on $h(X)$. Then:\\ [.05in] (1) A metric space
$(X,d)$ is hyperconvex if and only if every minimal extremal
function is a distance cone.\\(2) $(h(X),\rho)$ is well defined and
hyperconvex.\\(3) $X$ is isometrically embedded in $h(X)$, via the
map $d:X\rightarrow h(X)$ defined by $d_x(y)=d(x,y)$.\\(4) If
$X\subset A\subset h(X)$, then $h(A)$ is isometric to $h(X)$.\\(5)
If $f\in h(X)$ and the distance cone $h_v \in h(X)$, then
$\rho(h_v,f)=f(v)$.\\ (6) If we have $f\in h(X)$, then
$f(x)=sup_{w\in X} \{ d(x,w)-f(w)\}$.\\ (7) If $f \in h(X)$, then
$f$ is continuous. That is, we have $X\hookrightarrow h(X)
\hookrightarrow B(X)$ where the first mapping is the mapping $d$
defined in (3) the second map is the natural embedding of $h(X)$
into $B(X)$.
\end{Theorem}

\section {Metric Trees}

    In the following we denote the distance between two points $x,y \in X$ by $ xy:=d(x,y)$.
 \begin{Def}\label{definition}
 A \emph{metric tree} $X$ is a metric space $(X,d)$ satisfying the
following two axioms:\\[.05in]
(i) For every $x,y\in X$, $x\neq y$,  there is a uniquely determined
isometry  $$\varphi_{xy}:[0,d(x,y)]\rightarrow X$$ such that
$\varphi_{xy}(0)=x$, $ \varphi(d(x,y))=y $, and\\ [.05in](ii) For
every one-to-one continuous mapping $f:[0,1]\rightarrow X$ and every
$t \in [0,1]$, we have
$$
d(f(0),f(t))+d(f(t),f(1))=d(f(0),f(1)).
$$
\end{Def}

It is known \cite{ref6} that any metric tree $X$ has the four-point
property, but only a connected, complete metric space with the
four-point property is a metric tree. Since a metric tree is a space
in which there is only one path between two points $x$ and $y$, this
would imply that if $z$ is a point between $x$ and $y$ (that is, if
$xz+zy= xy $), then we know that $z$ is actually on the path between
$x$ and $y$. This motivates the next concept of a metric interval.\\

A \,\emph{metric interval} $<x,y>$ is defined as
$$<x,y>:=\{z \in X: xz+zy=xy\}.$$ Consider the function
$$h_x:<x,y>\rightarrow[0,xy]$$ defined by $h_x(z)=xz$. That is $h_x$
is the restriction of the distance cone to the metric interval\, It
was proved in \cite{ref5} that $(X,d)$ satisfying only \emph{the
first property of a metric tree} is equivalent to $h_x$ being a
\emph{bijective isometry}, which says that a metric interval is the
same as an interval in $\mathbb {R}$
\begin{Remark}
Suppose we have a metric segment $<x,y>$ in a metric tree. Since
metric trees satisfy the four-point property, if we take $u\in
<x,y>$\, and \,$u^*\in <x,y>$ we have $$ xu+uy=xy \,~\mbox{and}~
\,xu^*+u^*y=xy$$ third distance is $$uu^*+xy \leq xy$$ yielding
$uu^*=0$ or $u=u^*$. Thus the metric segment $<x,y> \subset
\{x,y,u\}$.
\end{Remark}

 We need the following three lemmas in order to prove Theorem 3.2
 below. Ideas behind these lemmas can be found in  \cite{ref5}. Nevertheless, we
 reconstruct and expand these ideas using Isbell's \cite{ref7} notation.
 Below in lemmas 3.2 and 3.3, we give a more detailed version of the
 proof given in \cite{ref5}. In \cite{ref5}, to prove the fact that
 the four-point property is inherited by the hyperconvex hull, the
 concept of ``thready spaces" was used which will be omitted in our discussion.

\begin{Lemma} (Dress)\\
(a) In a metric tree $(X,d)$, for any points $x$,$y$ \, and $z$\,
the intersection $$<x,y> \cap <x,z>$$ is a metric segment ending at
some point $u$.\\ [.05in] (b) In a metric tree $(X,d)$, we have
$$<x,y>\cap<y,z>\cap<z,x> \neq\emptyset$$ for all $x,y,z \in X$.
\end{Lemma}

Part (a) of the above lemma tells that if a portion of the metric
space looks like a line segment, and this segment splits into two,
the pieces can never connect again, so it must look rather like a
tree. Part (b) is expressing that metric trees are median.
\begin{Lemma} In a metric tree $(X,d)$ \, for $x,y \in X $ we have
$$<x,y>_X = <x,y>_{h(X)},$$ where $<x,y>_X= \{z\in X: \, xz+zy=xy\}.$
 Similarly $<x,y>_{h(X)}= \{z\in h(X): \, xz+zy=xy \}$.
\end{Lemma}
\textbf {Proof}: $X$, and therefore $h(X)$, are trees and it is
clear that $<x,y>_X \subset <x,y>_{h(X)}.$ To show the other
inclusion, consider the map $h_x:<x,y>_{h(X)}\rightarrow [0,xy]$\,
defined by\, $x\mapsto zx$. This is a bijective isometry since
$h(X)$ is a tree. On the other hand we also know that \, for all $ r
\in [0,xy]$, \, there exists $ x_r \in <x,y>_X $ \,  with  $ xx_r=r$
because $X$ is a tree. Therefore, if we take $z\in <x,y>_{h(X)}$, we
have $zx=wx$ for some $w \in <x,y>_X$. Therefore $h_x(z)=h_x(w).$
 Since $ h_x$ is injective, we have $z=w$. \hfill$\Box$
\begin{Lemma} If $(X,d)$ has the four-point property, then
$h(X)$ has the four-point property.
\end{Lemma}

\textbf {Proof}: First we show that if the metric space $(X,d)$ has
the four-point property, and if $f\in h(X)$, then $X \cup \{f\}$ has
the four-point property. Suppose $f,x,y,v \in X \cup \{f\}$. Then

\begin{eqnarray*}
xy+\rho (h_v,f) &=& xy+f(v) =\sup_{w\in X} \{
xy+vw-f(w)\} \\
&\leq& \max \left\{ \sup_{w \in X} \{ xv+yw-f(w) \} ,\, \sup_{w \in
X}
\{ xw+yv-f(w)\} \right\} \\
&=&\max \{ xv+f(y), \, yv+f(x) \}=\max\{xv+\rho (h_y,f),\, yv+\rho
(h_x,f) \}.
\end{eqnarray*}

This proves that $X \cup \{f\}$ has the four-point property. To
prove that $h(X)$ has the four-point property, use item (4) of
Theorem 2.3 and $X\subset X \cup \{f\} \subset h(X)$, which yields
$h(X \cup \{f\}) =h(X)$. Using the argument above, by taking $f_2
\in h(X\cup \{f_1\} )$, we see that $X\cup \{f_1,f_2\}$ has the
four-point property. Continuing in this manner and adding one point
at a time concludes the proof. \hfill $\Box$
\begin{Theorem}(Dress) A metric space is a metric tree if and
only if it is complete, connected and satisfies the four-point
property.
\end{Theorem}

\begin{Theorem}
Every complete metric tree is hyperconvex.
\end{Theorem}
\textbf{Proof}: Suppose $ (X,d)$ is a metric tree. Then by the above
theorem it has the four-point property, which in turn implies that
$h(X)$ has the four-point property. Since the hyperconvex hull is
connected, $h(X)$ is a metric tree as well. We would like to prove
that any minimal extremal function $f\in h(X)$ is a distance
cone.\,(i.e., $f$ has a zero).This is sufficient because, as
described in Remark 2.1, any pairwise overlapping collection can be
extended to a minimal extremal function, and this function having a
zero means that the point $x$ where $f(x)=0$ will be within the
radius of each closed ball in the original collection. In the
following we identify a point $x\in X$ with its isometric image
$h_x\in h(X)$. Start by fixing an $x\in X$ and use the minimality of
$f$ to obtain that, for each $\epsilon >0$, there is a point $y$,
depending on $\epsilon$, with $f(x)+f(y)\geq d(x,y)+\epsilon$.
Equivalently, for all $n\in N$, set $\epsilon=1/n$ and find $x_n$
with $f(x)+f(x_n)\geq d(x,x_n)+1/n$. Now, using Lemma 3.1 part (b.)
there is an element $g_n\in h(X)$ with
$$
 g_n\in <x,x_n>_{h(X)} \cap <x_n,f>_{h(X)} \cap
<f,x>_{h(X)}.
$$
This means $xg_n+g_nf=xf$ and $x_ng_n+g_nf=x_nf$, giving us
$$2g_nf+xg_n+g_nx_n=f(x)+f(x_n).$$ This equality is further reduced
to
$$2g_nf+xx_n=f(x)+f(x_n)$$ using the fact that $xg_n+g_n x_n=xx_n$.
Rewriting, we will have $$ g_nf=1/2(f(x)+f(x_n)-xx_n)\leq 1/2n .$$
We now use Lemma 3.2 to write $g_n \in <x,y>_{h(X)}=<x,y>_X$.
However, all elements of $<x,y>_X$ are distance cones, therefore
$g_n=h_{y_n}$ for some point $y_n \in X$, and $fg_n=f(y_n)$. Since
$fg_n\leq 1/2n$, we have a sequence of points $\{y_n\}$ with
$f(y_n)\leq 1/2n$. $\{y_n\}$ is a Cauchy sequence. Completeness
gives us a limit point $ y^\ast$ in $X$ and the continuity of $f$
implies $f(y^\ast)=0$. \hfill $\Box$
\begin{Remark}
There are two equivalent definitions of a metric tree. One
definition is due to A. Dress (named as T-theory).This definition
yields several ``properties" of metric intervals. The other
definition was given by J. Tits \cite{ref13} (named as $\mathbb
{R}$-trees), which  lists ``properties" of metric intervals as part
of the definition. W. A. Kirk \cite{ref10}, using J. Tits'
definition, proved that a metric space is a complete R-tree if and
only if it is hyperconvex and has unique metric segments. Here we
use A. Dress' definition  to show all complete metric trees are
hyperconvex. Moreover, Kirk's method of proof is quite different
then ours. Our aim is to use the elegant and geometrical nature of
\textbf{ the four-point property} for metric trees when making the
connection between hyperconvexity and metric trees.
\end{Remark}
\section {Extension Theorems and Metric Trees}
    The theory of Banach spaces could not have developed without the
  Hahn-Banach theorem. So it is natural to ask whether the same type of
 extension theorem is true in the context of metric spaces.
 This question have led Aronszajn  and Panitchpakti \cite{ref1}
 to the theory of hyperconvex spaces. They established the following
 theorem.
 \begin{Theorem}
 Let $X$ be a metric space. $X$ is hyperconvex if and only if every
 mapping $T$ of a metric space $Y$ into $X$ with some subadditive
 modulus of continuity $\delta(\epsilon)$ has, for any space $Z$
 containing $Y$ metrically, an extension $\tilde{T}:Z\rightarrow X$
 with the same modulus $\delta(\epsilon)$.
 \end{Theorem}
 It is worth noting that earlier L. Nachbin in \cite{ref12} proved a generalization of the Hahn-Banach
 theorem, stating that if the target space of a bounded linear map is
 an arbitrary real normed space, instead of the real numbers, then the
 extension is possible exactly when this target space is
 hyperconvex (he did not use the term ``hyperconvex").
 Extension theory for general bounded linear operators has a lot of
 unanswered questions even for basic cases. However, if one restricts
 the discussion to the extension of compact operators, there are a lot
 of elegant results (see \cite{ref13}).
 In the following, we discuss  $P_1$ spaces.

 \begin{Def}\label {definition}
 A metric space $(X,d)$ has the \emph{binary ball intersection
 property} if given any collection of closed balls that intersect
 pairwise, their total intersection is non-empty.
 \end{Def}
 It is clear that if a metric space is hyperconvex then it has the
 binary ball intersection property. For if the collection
 $\{B(x_i,r_i)\}$ intersects pairwise and if
 $x\in B(x_i,r_i)\cap B(x_j,r_j)$, then by the triangle inequality
$d(x_i,x_j)\leq d(x_i,x)+d(x,x_j)\leq r_i+r_j$ is satisfied. However
the binary ball intersection property
 does not imply hyperconvexity. If a space has the binary ball
 intersection property with the additional assumption that it is
 totally convex \cite{ref9}, then it is hyperconvex.
 \begin{Def}\label{definition}
 A Banach space $X$ is called \emph{1-injective}, or a $P_1$\emph{-space},
 if for every space $Y$ containing $X$ there is a projection $P$
 from $Y$ onto $X$ with $||P||\leq 1$.
 \end{Def}
 A real Banach space $X$ is $P_1$ if and only if it has the binary
 intersection property for balls, hence if and only if it is an
 absolute 1-Lipschitz retract  (see \cite{ref14}). The
 work of Nachbin, Goodner, Kelly and Hasumi characterizes real and
 complex $P_1$-spaces as the $C(K)$ spaces for extremally
 disconnected compact Hausdorff spaces $K$. For details see
 \cite{ref4}.

 An example of a $P_1$-space is a real $L_\infty(\mu)$ space with $\mu$ finite. This space has the binary
 intersection property, and hence it is a $P_1$-space. \\
 \begin{Theorem}
 Suppose X is a real Banach space that satisfies the four-point
 property. Then $X$ is a $P_1$-space.
 \end{Theorem}
 \textbf{Proof:} The result is a consequence of Theorem 3.1 and 3.2, and the fact
 that:
\begin{eqnarray*}
\lefteqn{\mbox{$X$ is hyperconvex}\Rightarrow }\\&& X\mbox{ has the
binary
 intersection property}\\&&\Leftrightarrow X \mbox{ is a
 $P_1$-space}\\
 &&\Leftrightarrow X\mbox{ is an absolute 1-Lipschitz
 retract.}
\end{eqnarray*}
\vspace{-.07in} \hfill$\Box$

 \begin{Remark}
  Matou\~{s}ek  in \cite{ref11} proves the following theorem. Let $Y$ be
 a metric tree and $X\subset Y$, and let $f$ be a mapping of $X$ into a
 Banach space $Z$ with Lipschitz constant $L$. Then, $f$ can be
 extended onto $Y$ with Lipschitz constant $CL$, where $C$ is an
 absolute constant. He also uses the four-point property in his
 proof.
 \end{Remark}


\begin{thebibliography}{99}
\bibitem{ref1} N. Aronszajn and P. Panitchpakdi, {\em Extension of
uniformly continous transformations and hyperconvex metric spaces},
Pacific J. Math. 6 (1956), 405--439.
\bibitem {ref2} D. Bugajewski and E. Grzelaczyk, {\em A fixed point
theorem in hyperconvex spaces}, Arch. Math. 75 (2000), 395--400.
\bibitem{ref3} P. Buneman, {\em A note on the metric properties of
trees}, J. Combin. Theory Ser. B, 17 (1974), 48--50.
\bibitem {ref4} M. M. Day, ``Normed Linear Spaces", Third edition,
Springer-Verlag, Berlin, Heidelberg, New York. 1973.
\bibitem{ref5} A. W. M. Dress, {\em Trees, tight extensions of metric
spaces, and the chomological dimension of certain groups: a note on
combinatorial properties of metric spaces}, Adv. in Math. 53 (1984),
321--402.
 \bibitem{ref6} A. W. M. Dress, V. Moulton and W. Terhalle,
{\em T-Theory, an overview}, European J. Combin. 17 (1996),
161--175.
\bibitem{ref7} J. R. Isbell, {\em Six theorems about injective metric
spaces}, Comment. Math. Helv. 39 (1964), 439--447.
\bibitem{ref8} W. B. Johnson, J. Lindenstrauss and D. Preiss, {\em Lipschitz
quotients from metric trees and from Banach spaces containing
${l}_1^1$ }, J. Funct. Anal. 194 (2002), 332--346.
\bibitem{ref9} M. A. Khamsi and W. A. Kirk,  ``An Introduction to
Metric Spaces and Fixed Point Theory", Pure and Applied Math.,
Wiley, New York, 2001.
\bibitem{ref10}W. A. Kirk, {\em Hyperconvexity of R-Trees}, Fund. Math.
156 (1998), 67--72.
\bibitem{ref11} J. Matou\v{s}ek, {\em Extension of Lipschitz
mappings on metric trees}, Comment. Math. Univ. Carolinae 31 (1990),
99--104.
\bibitem{ref12} L. Nachbin, {\em A theorem of  Hahn-Banach type
for linear transformations}, Trans. Amer. Math. Soc. 68 (1950),
28--46.
\bibitem{ref13} J. Tits,  {\em A theorem of Lie-Kolchin for trees},
Contributions to Algebra: a collection of papers dedicated to Ellis
Kolchin, Academic Press, New York, 1977.
\bibitem{ref14} M. Zippin, {\em Extension of bounded linear
operators}, Handbook of the geometry of Banach spaces, Vol.2 (2003),
1703--1741.
\end{thebibliography}
\end{document}